
\magnification=1200
\pretolerance=500 \tolerance=1000 \brokenpenalty=5000
\hoffset=0.5cm
\voffset=1cm
\hsize=12.5cm
\vsize=19cm
\parskip 3pt plus 1pt
\parindent=6.6mm


\font\fourteenbf=cmbx10 at 14.4pt
\font\twelvebf=cmbx10 at 12pt
\font\twelvebfit=cmmib10 at 12pt
\font\tenbfit=cmmib10
\font\sevenbfit=cmmib7
\font\eightrm=cmr8
\font\eightbf=cmbx8
\font\eighttt=cmtt8
\font\eightit=cmti8
\font\eightsl=cmsl8
\font\sevenrm=cmr7
\font\sevenbf=cmbx7
\font\sixrm=cmr6
\font\sixbf=cmbx6
\font\fiverm=cmr5
\font\fivebf=cmbx5

\font\tenCal=eusm10
\font\sevenCal=eusm7
\font\fiveCal=eusm5
\newfam\Calfam
  \textfont\Calfam=\tenCal
  \scriptfont\Calfam=\sevenCal
  \scriptscriptfont\Calfam=\fiveCal
\def\Cal{\fam\Calfam\tenCal}

\font\fourteenmsb=msbm10 at 14.4pt

\font\tenmsb=msbm10
\font\eightmsb=msbm8
\font\sevenmsb=msbm7
\font\fivemsb=msbm5
\newfam\msbfam
  \textfont\msbfam=\tenmsb
  \scriptfont\msbfam=\sevenmsb
  \scriptscriptfont\msbfam=\fivemsb
\def\Bbb{\fam\msbfam\tenmsb}


\font\eighti=cmmi8
\font\eightsy=cmsy8
\font\sixi=cmmi6
\font\sixsy=cmsy6

\skewchar\eighti='177 \skewchar\sixi='177
\skewchar\eightsy='60 \skewchar\sixsy='60

\catcode`\@=11

\def\pc#1#2|{{\bigf@ntpc #1\penalty\@MM\hskip\z@skip\smallf@ntpc #2}}

\def\tenpoint{%
  \textfont0=\tenrm \scriptfont0=\sevenrm \scriptscriptfont0=\fiverm
  \def\rm{\fam\z@\tenrm}%
  \textfont1=\teni \scriptfont1=\seveni \scriptscriptfont1=\fivei
  \def\oldstyle{\fam\@ne\teni}%
  \textfont2=\tensy \scriptfont2=\sevensy \scriptscriptfont2=\fivesy
  \textfont\itfam=\tenit
  \def\it{\fam\itfam\tenit}%
  \textfont\slfam=\tensl
  \def\sl{\fam\slfam\tensl}%
  \textfont\bffam=\tenbf \scriptfont\bffam=\sevenbf
  \scriptscriptfont\bffam=\fivebf
  \def\bf{\fam\bffam\tenbf}%
  \textfont\ttfam=\tentt
  \def\tt{\fam\ttfam\tentt}%
  \textfont\msbfam=\tenmsb
  \scriptfont\msbfam=\sevenmsb   \scriptscriptfont\msbfam=\fivemsb
  \def\Bbb{\fam\msbfam\tenmsb}%
  \abovedisplayskip=6pt plus 2pt minus 6pt
  \abovedisplayshortskip=0pt plus 3pt
  \belowdisplayskip=6pt plus 2pt minus 6pt
  \belowdisplayshortskip=7pt plus 3pt minus 4pt
  \smallskipamount=3pt plus 1pt minus 1pt
  \medskipamount=6pt plus 2pt minus 2pt
  \bigskipamount=12pt plus 4pt minus 4pt
  \normalbaselineskip=12pt
  \setbox\strutbox=\hbox{\vrule height8.5pt depth3.5pt width0pt}%
  \let\bigf@ntpc=\tenrm \let\smallf@ntpc=\sevenrm
  \normalbaselines\rm}
\def\eightpoint{%
  \let\yearstyle=\eightrm
  \textfont0=\eightrm \scriptfont0=\sixrm \scriptscriptfont0=\fiverm
  \def\rm{\fam\z@\eightrm}%
  \textfont1=\eighti \scriptfont1=\sixi \scriptscriptfont1=\fivei
  \def\oldstyle{\fam\@ne\eighti}%
  \textfont2=\eightsy \scriptfont2=\sixsy \scriptscriptfont2=\fivesy
  \textfont\itfam=\eightit
  \def\it{\fam\itfam\eightit}%
  \textfont\slfam=\eightsl
  \def\sl{\fam\slfam\eightsl}%
  \textfont\bffam=\eightbf \scriptfont\bffam=\sixbf
  \scriptscriptfont\bffam=\fivebf
  \def\bf{\fam\bffam\eightbf}%
  \textfont\ttfam=\eighttt
  \def\tt{\fam\ttfam\eighttt}%
  \textfont\msbfam=\eightmsb
  \scriptfont\msbfam=\sevenmsb   \scriptscriptfont\msbfam=\fivemsb
  \def\Bbb{\fam\msbfam\tenmsb}%
  \abovedisplayskip=9pt plus 2pt minus 6pt
  \abovedisplayshortskip=0pt plus 2pt
  \belowdisplayskip=9pt plus 2pt minus 6pt
  \belowdisplayshortskip=5pt plus 2pt minus 3pt
  \smallskipamount=2pt plus 1pt minus 1pt
  \medskipamount=4pt plus 2pt minus 1pt
  \bigskipamount=9pt plus 3pt minus 3pt
  \normalbaselineskip=9pt
  \setbox\strutbox=\hbox{\vrule height7pt depth2pt width0pt}%
  \let\bigf@ntpc=\eightrm \let\smallf@ntpc=\sixrm
  \normalbaselines\rm}


\newskip\LastSkip
\def\nobreakatskip{\relax\ifhmode\ifdim\lastskip>\z@
  \LastSkip\lastskip\unskip\nobreak\hskip\LastSkip
  \fi\fi}
\catcode`\;=\active
\catcode`\:=\active
\catcode`\!=\active
\catcode`\?=\active
\def;{\nobreakatskip\string;}
\def:{\nobreakatskip\string:}
\def!{\nobreakatskip\string!}
\def?{\nobreakatskip\string?}

\frenchspacing
\tenpoint

\newif\ifpagetitre                \pagetitretrue
\newtoks\hautpagetitre         \hautpagetitre={\hfil}
\newtoks\baspagetitre          \baspagetitre={\hfil}
\newtoks\auteurcourant       \auteurcourant={\hfil}
\newtoks\titrecourant         \titrecourant={\hfil}
\newtoks\baspagedroite       \newtoks\baspagegauche
\newtoks\hautpagegauche       \newtoks\hautpagedroite
\hautpagegauche={\tenbf\folio\hfill\eightrm\the\auteurcourant\hfill }
\hautpagedroite={\hfill\the\titrecourant\hfill\tenbf\folio }  
\baspagedroite={\hfil} \baspagegauche={\hfil}
\headline={\ifpagetitre\the \hautpagetitre\else\ifodd\pageno\the\hautpagedroite
\else\the\hautpagegauche\fi\fi}
\footline={\ifpagetitre\the\baspagetitre\global\pagetitrefalse
\else\ifodd\pageno\the\baspagedroite\else\the\baspagegauche\fi\fi}

\def\appeln@te{}
\def\vfootnote#1{\def\@parameter{#1}\insert\footins\bgroup\eightpoint
  \interlinepenalty\interfootnotelinepenalty
  \splittopskip\ht\strutbox 
  \splitmaxdepth\dp\strutbox \floatingpenalty\@MM
  \leftskip\z@skip \rightskip\z@skip
  \ifx\appeln@te\@parameter\indent \else{\noindent #1\ }\fi
  \footstrut\futurelet\next\fo@t}

\def\corners{\def\margehaute{\vbox to \hmargehaute{\hbox to \lpage
{\lefttopcorner\hfil\righttopcorner}\vfil}}
\def\margebasse{\vss\hbox to \lpage{\leftbotcorner\hfil\rightbotcorner}}}

\def\footnoterule{\kern-6\p@
  \hrule width 2truein \kern 5.6\p@} 

\catcode`\@=12

\def\pd#1#2 {\pc#1#2| }
\catcode`\@=11
\def\p@int{{\rm .}}
\def\p@intir{\discretionary{\rm .}{}{\rm .\kern.35em---\kern.7em}}
\def\pointir{\afterassignment\pointir@\let\next=}
\def\pointir@{\ifx\next\par\p@int\else\p@intir\fi\next}
\catcode`\@=12

\def\abstract#1{\vbox{\eightpoint \pc ABSTRACT|\pointir #1}}

\newdimen\srdim \srdim=\hsize
\newdimen\irdim \irdim=\hsize
\def\NOSECTREF#1{\noindent\hbox to \srdim{\null\dotfill ???(#1)}}
\def\SECTREF#1{\noindent\hbox to \srdim{\csname REF\romannumeral#1\endcsname}}
\newlinechar=`\^^J
\def\openauxfile{
  \immediate\openin1\jobname.aux
  \ifeof1
  \message{^^JCAUTION\string: you MUST run TeX a second time^^J}
  \let\sectref=\NOSECTREF
  \else
  \input \jobname.aux
  \message{^^JCAUTION\string: if the file has just been modified you may 
    have to run TeX twice^^J}
  \let\sectref=\SECTREF
  \fi
  \message{to get correct page numbers displayed in the table of contents^^J}
  \immediate\openout1=\jobname.aux}

\newcount\numbersection \numbersection=-1
\newcount\numberchapter \numberchapter=0
\def\markpage#1{
      \advance\numbersection by 1
      \immediate\write1{\string\def \string\CH
      \romannumeral\numberchapter
      #1\romannumeral\numbersection \string{%
      \number\pageno \string}}}

\def\titre#1|{{\null\baselineskip=18pt
                           {\fourteenbf
  \textfont\msbfam=\fourteenmsb
  \scriptfont\msbfam=\tenmsb
  \scriptscriptfont\msbfam=\sevenmsb
  \def\Bbb{\fam\msbfam\fourteenmsb}
                           \vskip 2.25ex plus 1ex minus .2ex
                           \leftskip=0pt plus \hsize
                           \rightskip=\leftskip
                           \parfillskip=0pt
                           \noindent #1
                           \par\vskip 2.3ex plus .2ex}}
                           \markpage{T}}

\def\section#1|{{\twelvebf
      \textfont0=\twelvebf  \scriptfont1=\tenbf  \scriptscriptfont1=\sevenbf
      \textfont1=\twelvebfit\scriptfont1=\tenbfit \scriptscriptfont1=\sevenbfit
      \par\penalty -500
      \vskip 3.25ex plus 1ex minus .2ex
      \noindent
      #1\vskip1.5ex plus 0.2ex}
      \markpage{S}}

\def\ssection#1|{{\bf
      \par\penalty -200
      \vskip 3.25ex plus 1ex minus .2ex
      #1\pointir}}

\def\skippage{\vfill\eject \ifodd\pageno \else \vphantom{}\vfill\eject \fi
              \numbersection=-1
              \advance\numberchapter by 1}

\long\def\th#1|#2\finth{\par\vskip5pt\noindent
              {\bf #1}{\sl \pointir #2}\par\vskip 5pt}

\def\remarque#1|{\par\vskip5pt\noindent{\bf #1}\pointir }
\def\rque#1|{\par\vskip5pt{\sl #1}\pointir}

\def\ieme{\raise 1ex\hbox{\pc{}i\`eme|}}
\def\omini{\raise 1ex\hbox{\pc{}o|}}
\def\emini{\raise 1ex\hbox{\pc{}e|}}
\def\ermini{\raise 1ex\hbox{\pc{}er|}}
\def\\{\mathop{\hbox{\tenmsb r}}\nolimits}

\newif\ifchiffre
\def\chiffre{\chiffretrue}
\chiffre
\newdimen\laenge
\def\lettre#1|{\setbox3=\hbox{#1}\laenge=\wd3\advance\laenge by 3mm
\chiffrefalse}
\def\article#1|#2|#3|#4|#5|#6|#7|%
    {{\ifchiffre\leftskip=7mm\noindent
     \hangindent=2mm\hangafter=1
\llap{[#1]\hskip1.35em}{\bf #2}\pointir {\sl #3}, {\rm #4}, \nobreak{\bf #5} ({\yearstyle #6}), \nobreak #7.\par\else\noindent
\advance\laenge by 4mm \hangindent=\laenge\advance\laenge by -4mm\hangafter=1
\rlap{[#1]}\hskip\laenge{\bf #2}\pointir
{\sl #3}, #4, {\bf #5} ({\yearstyle #6}), #7.\par\fi}}
\def\livre#1|#2|#3|#4|#5|%
    {{\ifchiffre\leftskip=7mm\noindent
    \hangindent=2mm\hangafter=1
\llap{[#1]\hskip1.35em}{\bf #2}\pointir{\sl #3}, #4, {\yearstyle #5}.\par
\else\noindent
\advance\laenge by 4mm \hangindent=\laenge\advance\laenge by -4mm
\hangafter=1
\rlap{[#1]}\hskip\laenge{\bf #2}\pointir
{\sl #3}, #4, {\yearstyle #5}.\par\fi}}

\def\divers#1|#2|#3|#4|%
    {{\ifchiffre\leftskip=7mm\noindent
    \hangindent=2mm\hangafter=1
     \llap{[#1]\hskip1.35em}{\bf #2}\pointir{\sl #3}, {\rm #4}.\par
\else\noindent
\advance\laenge by 4mm \hangindent=\laenge\advance\laenge by -4mm
\hangafter=1
\rlap{[#1]}\hskip\laenge{\bf #2}\pointir {\sl #3}, {\rm #4}.\par\fi}}
\def\dem{\par\noindent{\sl Proof}\pointir}

\def\carre{\hbox{\font\ppreu=cmsy10\ppreu\char'164\hskip-6.6666pt\char'165}}
\def\finpr{{\ \penalty 500\carre\par\vskip3pt}}


\def\frac#1/#2{\leavevmode\kern.1em
   \raise.5ex\hbox{\the\scriptfont0 #1}\kern-.1em
      /\kern-.15em\lower.25ex\hbox{\the\scriptfont0 #2}}

\def\buildo#1^#2{\mathrel{\mathop{\null#1}\limits^{#2}}}
\def\buildu#1_#2{\mathrel{\mathop{\null#1}\limits_{#2}}}

%
%
%
\def\aujour{\ifnum\day=1 1\ermini\else\number\day\fi\
\ifcase\month\or janvier\or f\'evrier\or mars\or avril\or mai\or juin\or
juillet\or aout\or septembre\or octobre\or novembre\or d\'ecembre\fi\
\number\year}
\def\today{\ifcase\month\or January \or February \or March \or April\or 
May\or June\or July\or August \or September\or October\or November\or
December\fi\ \number\day , \number\year}
\catcode`\@=11
\newcount\@tempcnta \newcount\@tempcntb 
\def\timeofday{{%
\@tempcnta=\time \divide\@tempcnta by 60 \@tempcntb=\@tempcnta
\multiply\@tempcntb by -60 \advance\@tempcntb by \time
\ifnum\@tempcntb > 9 \number\@tempcnta:\number\@tempcntb
  \else\number\@tempcnta:0\number\@tempcntb\fi}}
\catcode`\@=12

\def\Gr{\mathop{\rm Gr}\nolimits}

\def\div{\mathop{\rm div}\nolimits}
\def\mod{\mathop{\rm mod}\nolimits}

\def\lra{\longrightarrow}
\def\ld{,\ldots,}

\def\wt{\widetilde}

\def\bC{{\Bbb C}}
\def\bR{{\Bbb R}}
\def\bP{{\Bbb P}}
\def\bN{{\Bbb N}}
\def\bG{{\Bbb G}}
\def\bZ{{\Bbb Z}}
\def\bQ{{\Bbb Q}}

\def\cF{{\Cal F}}

\def\cO{{\Cal O}}

\def\rk{\mathop{\rm rank}\nolimits}

\def\Pic{\mathop{\rm Pic}\nolimits}

\def\longhook{\lhook\joinrel\longrightarrow}

\titre{Logarithmic Jets and Hyperbolicity}|
\auteurcourant{Jawher El Goul}
\titrecourant{\eightpoint Logarithmic Jets and Hyperbolicity  }

\centerline{by  Jawher El Goul 
(Toulouse III)}

\vskip10pt
June 1st, 2000, printed on \today, \timeofday

\vskip 0.5cm
{\eightpoint\baselineskip=10pt{
\noindent
{\bf Abstract.}
We prove that  the complement of a very generic curve of degree $d$ at least equal
to $15$ in ${\Bbb P}^2$ is hyperbolic in the sens of Kobayashi (here, the 
terminology ``very generic'' refers to complements of countable unions of proper 
algebraic subsets of the parameter space). We first consider the
Dethloff and Lu's generalisation  to the logarithmic situation of
Demailly's jet bundles. We study their base loci for surfaces of 
log-general type in the same way as it was done in the compact case by 
Demailly and El Goul. With some condition on log-Chern classes, any entire 
holomorphic map to the surface can be lifted as a leaf of some foliation on 
a ramified covering.  
Then we obtain a logarithmic analogue of   
McQuillan's result on holomorphic foliations which permits to conclude. 
Using the logarithmic formalism, we even obtain some simplifications of the 
original proof in McQuillan's work.
}}
\vskip 0.2cm

\section 0. Introduction|

In 1970 S. Kobayashi [Ko70] posed the following problems: Is it true 
that the complement of generic hypersurface of degree $d\geq e(n)$ in 
${\Bbb P}_\bC^{n}$ is hyperbolic? Is this true for $e(n)=2n+1$?   
Later Green [Gr75] for $n=2$ and Zaidenberg [Za87] for arbitrary $n$ 
proved that for $d\leq 2n$ those complements contain $\bC^\star$ and then
are not hyperbolic.

In this paper we will study the case of complements of smooth
curves in ${\Bbb P}_\bC^{2}$ where for $d\geq 4$ this is equivalent
to the nonexistence of nonconstant entire curves by Brody Reparametrisation 
Lemma [Br78]. When the curve is supposed to have many 
components this had been studied by many authors, see [DSW92,94] for a 
complete bibliography and the study of the case of three components
(see also the recent work [BD00]). 
In the smooth case, after an example of hyperbolic complement given 
by Nadel in [Na89] (for $d\ge
21$), Zaidenberg [Za89] gave examples for all $d\geq 5$.

The first positive answer to this question was given in the work of
Siu and Yeung [SY95], the bound they obtain is quite high. Their method 
consists of an explicit construction of special second order differential
operators on an associated surface in ${\Bbb P}_\bC^{3}$  ramified over ${\Bbb
  P}_\bC^{2}$. This was done by an imitation of the construction of holomorphic
1-forms on Riemann surfaces and a clever reduction of the problem to a 
resolution of linear systems. Those operators are such that their pullbacks 
by the lifting of every entire curve must vanishes identically. This follows 
immediately from an Ahlfors type result.

In [DEG00], after studying the compact analogue of the above  
conjecture and proving that a generic suface of degree $d\geq 21$ in
${\Bbb P}_\bC^{3}$
is hyperbolic, Demailly and the author, using the same covering trick,
obtained the bound $21$ also for complementary of curves in ${\Bbb P}_\bC^{2}$.
This is made by using the whole force of Demailly's jet bundles introduced in
[De95] and a McQuillan's result on holomorphic foliations [Mc98].

Here we obtain the bound $15$. Our main theorem is the following

\th Main Theorem|
The complement of a very generic curve in ${\Bbb P}_\bC^{2}$ is hyperbolic
and hyperbolically imbedded for all degrees $d\geq 15$.
\finth

We follow almost the same strategy as in [DEG00] 
with the difference that we use an logarithmic analogous package introduced 
in [DL96], which introduces some additionnal technical complications.
Dethloff and Lu's jet bundles are a compactification \`a la Demailly of
Noguchi's logarithmic jet bundles introduced in [No84]. Using Riemann-Roch
and a refined study of the base loci associated to those jet bundles,
we reduce the problem to the study of holomorphic foliations on log-general
type surfaces. We prove that such foliations do not admit a parabolic
Zariski-dense leaf. We generalize, in particular, McQuillan's result [Mc98] 
on Green-Griffiths conjecture to 
log-general type surfaces with $\overline{c}_1^2 >\overline{c}_2.$ 
This logarithmic point of view permits the observation that 
McQuillan's refined tautological inequlity is actually an easy consequence
of a logarithmic tautological inequality obtained by [Vo99] and proved in the
same way as the simple one  (see also [Br00]). 
 
The paper is organized as follows: In section 1 we first recall the main 
definitions and results in [DL96]. Then we introduce the $2$-jet threshold 
 of a log-general type surface $(X,C)$. We consider the case when the  
Picard group is $\bZ$ and we construct, with some conditions on log-Chern 
classes and  the $2$-jet threshold, a ramified cover $\tilde{X}$ on  $X$ 
such that every entire curve in the complement $X\setminus C$ 
could be lift as a leaf of a holomorphic foliation on $\tilde{X}$.
The last part of this section is devoted to estimate the $2$-jet threshold
in the case of $({\Bbb P}_\bC^{2},C)$ where $C$ a generic smooth plane curve.

In section 2 we study foliations on log-general type surfaces as in [Mc98].
The method we adopt is parallel to Brunella's work [Br99]. We obtain
that those foliations do not have a Zariski dense parbolic leaf.

The author would like to thank warmly Professors Marco Brunella, Jean-Pierre 
Demailly and Gerd Dethloff for interesting remarks and suggestions.

\section 1. Logarithmic Demailly jet bundles|

\ssection 1.1. Background material| Here we will consider a logarithmic 
generalisation of Demailly's invariant jets introduced in [De95], 
this is done by \hbox{ Dethloff and  Lu} in [DL96]. 

Let $X$ be an $n$-dimensional complex manifold with a normal crossing divisor
$D$. According to  \hbox{Iitaka} [Ii77], the logarithmic cotangent sheaf
$\overline{T}^\star_X=T^\star_X{(\log{D})}$ is defined to be the locally free sheaf generated by
$T^\star_X$ and the logarithmic differentials ${ds_j}/{s_j}$ , where $s_j =0$ 
is a local equation for the
irreducible components of $D$. Its dual the logarithmic tangent sheaf 
is the sheaf of germs of 
vector fields tangent to $D$, denoted by $\overline{T}_X= T_X {(-\log{D})}$. 

Recall from [GG80] that the $k$-jet bundle $J_kX$ is defined as the set 
of equivalence classes of holomorphic maps $f:(\bC,0)\to(X,x)$, with the 
equivalence relation $f\sim g$ if and only if  they have the same Taylor
expansions of order $k$ in some local coordinate system of $X$ near~$x.$
We denote the equivalence class of $f$ by $j_k (f)$.
In [No84], Noguchi generalised this object to the logarithmic situation 
as follows. Let $\omega \in H^0(U, T^\star_X )$ be a holomorphic section over 
an open subset $U\subset X$. For a germ of a holomorphic map $f$ in $U$ we put
$f^\star\omega = Z(t)dt$. Then we have a well defined holomorphic mapping
 $$\tilde\omega : J_kX|_U \to \bC^k\ ; \ j_k(f)\to (Z^{(j)}(0))_{0\leq j\leq k-1}.$$
Now we say that a holomorphic section $s \in H^0(U,J_k X)$ is a logarithmic
$k$-jet field if the map $\tilde\omega\circ s |_V : V \to \bC^k$ is
holomorphic for all $\omega \in H^0(V, \overline{T}^\star_X )$ and for all
open subset $V$ of $U$. The set of logarithmic $k$-jet fields over open
subsets of $X$ defines a subsheaf of $J_k X$ called the logarithmic $k$-jet
bundle of $(X,D)$, which we denote by $\overline{J}_k X$.

In [DL96], Dethloff and Lu constructed a more geometrically relevant
jet bundles (in the same way as done in [De95] for the non-logarithmic case)
by considering a suitible "quotient" of this bundle
by the action of the group $\bG_k$ containing all  
germs of $k$-jets biholomorphisms of $(\bC,0)$, that is, the group of
germs of biholomorphic maps
$$
t\mapsto\varphi(t)=a_1t+a_2t^2+\cdots+a_kt^k,\qquad
a_1\in\bC^\star,~a_j\in\bC,~j\ge 2.
$$
As a generalization of Demailly's directed jets to the logarithmic context,
\hbox{Dethloff and Lu} defined a log-directed manifold to be the
triple $(X,D,V)$ where $V$ is a holomorphic  subbundle of $\overline{T}_X
$ of rank $r$. 
To the log-directed manifold $(X,D,V)$, one associates inductively a sequence of
directed manifolds $(\overline{X}_k,D_k,V_k)$ as follows. Starting with $(\overline{X}_0,D_0,V_0) = (X,D,V)$, 
one puts inductively $ \overline{X}_k= P(V_{k-1})$ with its natural projection $\pi_k$
to $\overline{X}_{k-1}$ (where $P(V)$ stands for the projectivized bundle of lines in the vector 
bundle $V$), where $D_k= \pi_k^\star(D_{k-1})$ and $ V_k$ is the subbundle of 
${T_{ \overline{X}_k}}(-\log{D_k})$ defined at any point $(x,[v])\in \overline{X}_k$, $v\in V_{k-1,x}$, by
$$ V_{k,(x,[v])} = \Big\{ \xi \in {T_{{ \overline{X}_k},{(x,[v])}}(-\log{D_k}) }~;~ {(\pi_k)}_\star\xi
\in \bC\cdot v\Big\},~$$
with $ \bC\cdot v\subset V_{k-1,x} \subset T_{\overline{X}_{k-1},x}(-\log{D_{k-1}})~.$
                                                               
We denote by $\cO_{\overline{X}_k}(-1)$ the tautological line subbundle of
$\pi^\star_k V_{k-1}$, such that 
$$\cO_{\overline{X}_k}(-1)_{ (x,[v])} = \bC \cdot v,$$
for all $(x,[v])\in \overline{X}_k = P(V_{k-1})$. By definition, the bundle $V_k$ fits
in an exact sequence
$$
0\lra T_{\overline{X}_k/\overline{X}_{k-1}}\lra V_k\buildo{\lra}^{\pi_{k\star}}
\cO_{\overline{X}_k}(-1)\lra 0,
$$
and the Euler exact sequence of $T_{\overline{X}_k/\overline{X}_{k-1}}$ yields
$$
0\lra\cO_{\overline{X}_k}\lra \pi_k^\star V_{k-1}\otimes\cO_{\overline{X}_k}(1)
\lra T_{\overline{X}_k/\overline{X}_{k-1}}\lra 0.
$$

From these sequences, we infer
$$\rk V_k = \rk V_{k-1}=\cdots=\rk V=r,\qquad
\dim \overline{X}_k = n+k(r-1).$$
We note
$$
\pi_{k,j}=\pi_{j+1}\circ\cdots\circ\pi_{k-1}\circ\pi_k:\overline{X}_k \lra \overline{X}_j,
$$
be the natural projection.

The canonical injection $\cO_{\overline{X}_k}(-1) \hookrightarrow \pi^\star_k
V_{k-1}$ and the exact sequence
$$
0\lra T_{\overline{X}_{k-1}/\overline{X}_{k-2}} \lra V_{k-1} \buildo \lra^{
{{(\pi_{k-1})}_\star}}
\cO_{\overline{X}_{k-1}}(-1) \lra 0
$$
yield a canonical line bundle morphism
 $$\cO_{\overline{X}_k}(-1) \buildo \longhook^{\ {(\pi_{k,k-2})}^\star \circ
{(\pi_{k-1})}_\star } \pi^\star_k\  \cO_{\overline{X}_{k-1}}(-1)$$
which admits precisely the hyperplane section $ \Gamma_k := P(T_{\overline{X}_{k-1}/\overline{X}_{k-2}})\subset
\overline{X}_k=P(V_{k-1})$  as its zero divisor. Hence we find
$\cO_{\overline{X}_{k}} (-1) = \pi^\star_{k}\  \cO_{\overline{X}_{k-1}}(-1) \otimes
\cO(-\Gamma_k)$ and using the notation $\cO_{\overline{X}_k}(a_1,a_2):=\pi_{k}^\star
\cO_{\overline{X}_{k-1}}(a_1)\otimes\cO_{\overline{X}_k}(a_2)$, 
$$
\cO_{\overline{X}_k}(-1,1)\simeq\cO (\Gamma_k)
$$ 
is associated with an effective divisor in~$\overline{X}_k$. 

For simplicity let us consider the case $V=\overline{T}_X$
and let \hbox{$f:\Delta_r\to X\setminus D $} tangent to $V$ be a nonconstant 
trajectory. Then $f$ lifts to a well defined and unique trajectory 
$f_{[k]}:\Delta_r\to \overline{X}_k\setminus D_k $ of $\overline{X}_k$ tangent 
to $V_k$.
Moreover, the derivative $f_{[k-1]}'$ gives rise to a section
$$f_{[k-1]}':T_{\Delta_r}\to f_{[k]}^\star\cO_{\overline{X}_k}(-1).$$
With any section $\sigma$ of $\cO_{\overline{X}_k}(m)$, $m\ge 0$, on any open
set $\pi_{k,0}^{-1}(U)$, $U\subset X\setminus D$, we can associate a 
holomorphic differential operator $Q$ of order $k$ acting on $k$-jets of 
germs of curves
$f:(\bC,0)\to U$ tangent to $V$, by putting
$$
Q(f)(t)=\sigma(f_{[k]}(t))\cdot f_{[k-1]}'(t)^{\otimes m}\in\bC.
$$
From [De95] this correspondence is, in fact, bijective. To see
what happen with logarithmic jets recall the following characterisation of
log-jet differentials in [DL96]: 
\th 1.1.1. Proposition ([DL96]) |
A holomorphic function $Q$ on $\overline{J}_k X|_U$ on some connected
open subset $U\subset X$ which satisfies
 $$(*)\quad Q(j_k(f\circ \phi))=\phi'(0)^m Q(j_k (f)) \quad\forall  j_k (f) \in J_kX|_V
 \quad\rm{and}\quad \forall \phi\in\bG_k$$
over some open subset $V$ of $U\setminus D$ defines a holomorphic
section of $\cO_{\overline{X}_k}(m)$ over $U,$ and vice versa.
\finth

Now, using the characteristion above, the definition of $\overline{J}_k X$
(in fact, on the component $D_j$ of the divisor $D$ the $k$-th derivative of 
the function $\log s_j(f)$ is holomorphic on $\overline{J}_k X$) and the fact 
that holomorphic functions satisfying $(*)$ for all
$\phi\in \bC^\star\subset \bG_k$ are homogenous polynomials, we obtain

\th 1.1.2. Proposition| 
The direct image $(\pi_{k,0})_\star \cO_{\overline{X}_k}(m)$ coincides with
the sheaf $\cO(E_{k,m}{\overline{T}^\star}_X)$ of logarithmic jet 
differentials, that is, the locally free sheaf generated by all polynomial
operators in the derivatives of order $1,2,\cdots, k$ of $f$, together
with the extra function $\log {s_j}(f)$ along the $j$-th component of $D$,
which are moreover invariant under arbitrary
changes of parametrization: a germ of operator
$Q\in E_{k,m}{\overline{T}^\star}_X $ is 
characterized by the condition that, for every germ $f$ in $X\setminus D$ and
every germ $\varphi$ of $k$-jet biolomorphisms of  $(\bC,0)$,
$$
Q(f\circ\varphi)=\varphi^{\prime m}\;Q(f)\circ \varphi.
$$
\finth

A basic result from [DL96] relying on the Ahlfors-Schwarz lemma, is the 
following, for the $1$-jet case see [No77] and [Lu91].

\th 1.1.3. Theorem ([DL96])|
 If $(X,D)$ has a $k$-jet metric $h_k$, i.e. a singular metric in the sens of
Demailly on $\cO_{\overline{X}_k}(-1)$ with negative curvature (along $V_k$),
then every entire curve  $f:\bC\to X\setminus D$ is such that $f_{[k]}(\bC)
\subset \Sigma_{h_k}$, where $\Sigma_{h_k}$ denotes the singular set of $h_k$.
\finth

An important case where the previous theorem applies is when  there are
some  integers $k,m>0$ and an ample line bundle $A$ on $X$ such that 
$$H^0(\overline{X}_k,\cO_{\overline{X}_k}(m)\otimes
{(\pi_{k,0})}^\star A^{-1})\simeq H^0(X,E_{k,m}{\overline{T}^\star}_X \otimes A^{-1})$$
has nonzero sections $\sigma_1\ld\sigma_N$. Then, we can construct 
$k$-metric of negative curvature, singular on their base locus 
$Z\subset \overline{X}_k$. 

By definition, a line bundle $L$ is {\it big} if there exists an ample divisor 
$A$ on $X$ such that $L^{\otimes m}\otimes \cO(-A)$ admits a nontrivial
global section when $m$ is large (then there are lot of sections, namely
$h^0(X,L^{\otimes m}\otimes \cO(-A))\gg m^n$ with $n=\dim X$).

As a consequence, Theorem 1.1.3 can be applied when $\cO_{\overline{X}_k}(1)$
is big or its restriction on some subvariety is big 
(See Theorem 4.3 in [DL96]).

In view of studying degeneration of entire curve
drawn on a variety of log-general type and of Theorem~1.1.3, it is especially
interesting to compute the base locus of the global sections of 
log-jet differentials, that is, the intersection 
$$
\overline{B}_k:=\bigcap_{m>0} \overline{B}_{k,m} \subset \overline{X}_k
$$
of the base loci $\overline{B}_{k,m}$ of all line bundles
$\cO_{\overline{X}_k}(m)\otimes\pi^\star_{k,0} \cO (-A)$, where $A$ is a given
arbitrary ample divisor over $X$.
 
\remarque 1.1.4. Remark| The $1$-jet case was studied in [Lu91]. 
Using Riemann-Roch [Hi66], we prove that 
if $(X,D)$ is a nonsingular surface of log-general type with logarithmic 
Chern classes $\overline{c}_1^2>\overline{c}_2$ then there is  a lot of 
log-symmetric differentials, i.e. sections in $E_{1,m}\overline{T}^\star =
S^m{\overline{T}^\star_X}$, and the base locus $\overline{B}_1$ is 
2-dimentionnal. Unfortunately, the ``order 1'' techniques  
are insufficient to deal with complement of smooth curve $C$ of degree $d$
in ${\Bbb P}^2$,  because in this case $$\overline{c}_1^2=(d-3)^2<\overline{c}_2=(d^2-3d+3).$$
Lemma 1.4.1 below shows in fact that $H^0(X, S^m {\overline{T}^\star_X})=0$ for all $m>0$.

\ssection 1.2. Base locus of logarithmic 2-jets|From now on, we suppose that 
$(X,D)$ is a nonsingular minimal surface of log-general
type (i.e. with $\overline{K}_X:=  K_X\otimes\cO(D)$ big and nef)  
and let us study the base locus $ \overline{B}_2 $ in $\overline{X}_2$.
As in the non-logarithmic case the bundle of log-jet differentials of order
$2$ has the following filtration
$$ 
\Gr^\bullet E_{2,m}{\overline{T}^\star_X}=\bigoplus_{0\le j\le m/3}
S^{m-3j}{\overline{T}^\star_X}\otimes\overline{K}_X^{\otimes j}.
$$
This filtration consists in writing an invariant polynomial log-differential 
operator outside $D$ as
$$
Q(f)=\sum_{0\le j\le m/3}\quad\sum_{\alpha\in\bN^2,\,|\alpha|=m-3j}
a_{\alpha,j}(f)\,(f')^\alpha(f'\wedge f'')^j
$$
where 
$$
f=(f_1,f_2),\qquad (f')^\alpha=(f'_1)^{\alpha_1}(f'_2)^{\alpha_2},
\qquad f'\wedge f''=f'_1f''_2-f''_1f'_2.
$$
On a component $D_j$ of $D$ given in local coordinate by $z_1=0$ we replace
only $f_1$ by $\log{f_1}$ in this expression.
A~calculation based on the above filtration of $E_{2,m}\overline{T}^\star_X$ and
 Riemann-Roch yields
$$\chi\big(X,E_{2,m}\overline{T}^\star_X\big)={m^4 \over {648}}(13\overline{c}_1^2-9\overline{c}_2) + 
O(m^3 ).$$   
On the other hand,
$$
H^2(X,E_{2,m}\overline{T}^\star_X\otimes\cO (-A))=
H^0(X,K_X\otimes E_{2,m}\overline{T}_X\otimes\cO(A))
$$
by Serre duality. From the filtration above, $K_X\otimes (E_{2,m}\overline{T}_X)\otimes\cO(A)$
admits a filtration with graded pieces
$$
S^{m-3j}\overline{T}_X\otimes\overline{K}_X^{\otimes -j}\otimes K_X\otimes\cO(A).
$$
Recall now that $\overline{T}_X$ is semi-stable (see [KR85] and [TY87]) so by Bogomolov's 
vanishing theorem [Bo79], we have $h^0 \big( X, S^p \overline{T}_X\otimes
\overline{K}_X^{\otimes q} \big)=0$, $p-2q>0$. This implies that
$$h^2 \big( X,E_{2,m}\overline{T}^\star_X\otimes\cO(-A)) =0$$
for $m$ large. Consequently we get the following 

\th 1.2.1. Theorem |
If $(X,D)$ is an algebraic surface of general type and
$A$ an ample line bundle over~$X$, then
$$
h^0(X,E_{2,m}\overline{T}^\star_X\otimes\cO(-A))\ge
{m^4\over 648}(13\,\overline{c}_1^2-9\,\overline{c}_2 )-O(m^3).
$$
In particular$\,:$
 If $13\,\overline{c}_1^2-9\,\overline{c}_2>0$, then $\overline{B}_2\neq \overline{X}_2$.
 
\finth

In the special case when $X={\Bbb P}^2$ and $D=C$ is a smooth plane curve
of degree $d$, we take $A= \cO_{{\Bbb P}^2}(1)$. Then we have 
$\overline{c}_1=(3-d)h$ and $\overline{c}_2=(d^2-3d+3)h^2$ where $h=c_1(\cO_{{\Bbb P}^2}(1))$,
$h^2=1$, thus  
$$\chi\big( E_{2,m}\overline{T}^\star_X \otimes \cO (-A) \big)= 
(4\,d^2-51\,d+90){m^4 \over {648}}+ O(m^3 ).$$
A straightforward computation shows that the leading coefficient
$4\,d^2-51\,d+90$ is positive if $d\geq 11$. Thus, we obtain

\th 1.2.2. Corollary|
 For every smooth curve  of degree $d\ge 11$  the associated log-surface
$(\bP^2, C)$ has its $2$-jets base locus  $\overline{B}_2\neq \overline{X}_2$.
\finth

\remarque 1.2.3. Remark|As a consequence of the calculus above and with the 
previous condition on Chern classes, every holomorphic entire curve $f$ into  
$X\setminus D$ could be lifted in $\overline{X}_2$ and its image is contained in an 
irreducible component $Z$ of $\overline{B}_2$. We have to distinguish three cases
\smallskip
\item{\rm(a)} $\pi_{2,1}(Z)=\overline{X}_1$, then $Z$ is three dimensional and called
 horizontal, in this case the $2$-jet lifting of $f$ is a leaf of a foliation
 by curves on $Z$. In fact the lifting of $f$ is tangent to $T_Z\cap V_2$
 which defines a distribution of lines on a Zariski open subset
 of $Z$ which is obviously integrable. 
\smallskip
\item{\rm(b)} $\pi_{2,0}(Z)= X$ and we are not in case (a), then our curve
 $f$ could be lifted to $\overline{X}_1$ as a leaf of foliation by curves on 
 the surface
 $Y:=\pi_{2,1}(Z)$ (defined by the distribution $T_Y\cap V_1$)
 
\smallskip
\item{\rm(c)} The curve $f$ is degenerated.
\medskip\noindent 
 
The difficulty, in the case (a), is to  study (singular) foliation by
curves on a variety of dimension bigger than two. Actually, we have no
reasonable model with reducible foliated singularities in this case
until now (see however [Mc99]). So the next step is to show that
with a slightly stronger condition on Chern classes, in the case of 
$(\bP^2, C)$, we have to consider only foliations on surfaces.

\ssection 1.3. Existence of the multi-foliation|Our aim now is to study the 
restriction of the tautological line bundle
on 2-jets on a $3$-dimentionnal horizontal component $Z$ of $\overline{B}_2$.
Let us first make the following useful definition (as in [DEG00]) 

\th 1.3.1 Definition |
Let $(X,D)$ be a nonsingular projective variety of log-general type. 
We define the {\rm $k$-jet log-threshold} $\overline\theta_k$ of $(X,D)$ to be the 
infimum
$$
\overline\theta_k=\inf_{m>0}\overline\theta_{k,m} \in\bR,
$$ 
where $\overline\theta_{k,m}$ is the smallest rational number $t/m$ such that there 
is a non zero section in $H^0(X,E_{k,m}\overline{T}^\star_X\otimes\cO(t\,\overline{K}_X))$
$($assuming that $t\,\overline{K}_X$ is an integral divisor, $t\in\bQ)$.
\finth

In the case when the Picard group equals $\bZ$ we have a more clear
idea about the jet log-threshold. In fact we have

\th 1.3.2. Lemma|
Let $(X,D)$ be a nonsingular surface of log-general type with  $\Pic(X)=\bZ.$ 
Suppose that $\overline\theta_1\geq 0$ and $\overline\theta_2<0$. 
Then $$ \overline{B}_2 \subset Z_\sigma = Z\cup\Gamma_2,$$ 
where $Z$ is an irreducible section and $ Z_\sigma $ is the set of
zeros of a section $\sigma$
in
$H^0(\overline{X}_2,\cO_{\overline{X}_2}(m_0)\otimes\cO(t_0\,\overline{K}_X)).$
Moreover, in the case $\overline{B}_2 = Z_\sigma$, we have $\overline\theta_2={t_0/m_0}$.

\finth
\dem
As $\overline\theta_2<0$ we have a nontrivial section 
$$
s\in H^0(X,E_{2,m}\overline{T}^\star_X\otimes\cO(t\,\overline{K}_X)),
\qquad m>0,~~t\in\bQ_-
$$
Let $u_1=\pi_{2,1}^\star\cO_{\overline{X}_1}(1)$ and
$u_2=\cO_{\overline{X}_2}(1),$
then its zero divisor
$$
Z_s=m\,u_2 + t\,\pi_{2,0}^\star K_X\qquad\hbox{in $\Pic(X_2)$},
$$
Let $Z_s=\sum p_jZ_j$ be the decomposition of $Z_s$ in 
irreducible components. From the equality 
$\Pic(\overline{X}_2)=\Pic(X)\oplus\bZ u_1\oplus \bZ u_2$ and the 
assumption $\Pic(X)\simeq\bZ$, we find 
$$Z_j\sim a_{1,j} u_1 + a_{2,j} u_2 + t_j\,\pi_{2,0}^\star \overline{K}_X,$$
for suitable integers $a_{1,j},\ a_{2,j} \in \bZ$ and rational 
numbers $t_j\in \bQ$. We can prove that (see Lemma 3.3 in [DEG00]), as $Z_j$ 
is effective,
we must have one of the following three disjoint cases:
\smallskip
\item{$\scriptstyle\bullet$} $(a_{1,j},a_{2,j})=(0,0)$ and 
$Z_j\in\pi_{2,0}^\star\Pic(X)$, $t_j>0\,$;
\smallskip
\item{$\scriptstyle\bullet$} $(a_{1,j},a_{2,j})=(-1,1)$, then $Z_j$ contains
$\Gamma_2$, so $Z_j=\Gamma_2$ and $t_j=0\,$;
\smallskip
\item{$\scriptstyle\bullet$} $a_{1,j}\ge 2 a_{2,j}\ge 0$ and 
$m_j:=a_{1,j}+a_{2,j}>0$. 
\smallskip
We can suppose $t_0/m_0 =\min{t_j\over m_j}$ then $t_0$ is clearly negative.
Now we have $a_{2,0}\neq 0 $ because $\theta_1\geq 0$ and  then $Z_0$ gives 
a section 
$$\sigma\in 
H^0\big(\overline{X}_2,\cO_{\overline{X}_2}(m_0)\otimes \pi_{2,0}^\star
\cO(t_0\overline{K}_X)\big)$$
(we use the identity $\cO_{\overline{X}_2}(a_1,a_2)=\cO_{\overline{X}_2}(a_1+a_2)\otimes\cO_{\overline{X}_2}(-a_1\Gamma_2)$).
Then, by definition, we obtain $\overline{B}_2 \subset Z_\sigma$ and 
we have equality if and only if $Z$ is the unique irreducible section
with $t<0$. As $t_j/m_j\leq t/m$ we conclude, in this case that
$\overline\theta_2={t_0/m_0}.$
\finpr

As a generalisation to the log-case of the main theorem of [DEG00] we have
the following

\th 1.3.3. Theorem|
Let $(X,D)$ be a nonsingular surface of log-general type with  $\Pic(X)=\bZ .$ 
Suppose that $\overline\theta_2<0$ and that the log-Chern numbers of $X$ 
satisfy~
$$\smash{\displaystyle
{ (13+12\,\overline\theta_2)\overline{c}_1^2 > 9\,\overline{c}_2.}}
$$
Then every Zariski-dense holomorphic map $f : \bC \to X\setminus D$ is a
leaf of an algebraic multi-foliation on $X$.

\finth

\dem 
If $\overline\theta_1<0$ then $\overline{B}_1\neq \overline{X}_1$ and we 
conclude by a direct application of Theorem 1.1.3.(the foliation is defined by
 the intersection of $V_1$ and the tangent of the irreducible component 
 of $\overline{B}_1$ which
 contains the lifting of $f$ to $1$-jets) so we suppose that
$\overline\theta_1\geq 0$. 
As $\overline\theta_2<0$, then $\overline{B}_2 \neq \overline{X}_2$
and the discussion made in Remark 1.2.3. shows that
we have to consider only the case when the lifting of $f$ to $2$-jets is 
contained in
a horizontal irreducible divisor $Z$ in $\overline{X}_2$. 
By Lemma 1.3.2 we have
$$
Z\sim a_1 u_1+a_2 u_2 +t_0 \pi_{2,0}^\star \overline{K}_X\qquad\hbox{in $\Pic(\overline{X}_2)$},~~t_0\in\bQ,~~t_0<0,
\qquad a_1+a_2=m_0,$$ 
where ${t_0/ m_0}=\overline\theta_2.$
 
Our aim now is to prove that the restriction of the tautological line bundle
to $Z$ is big. First, we have the following intersection equalities
$$
\eqalign{
&u_1^4=0,\quad u_1^3u_2=\overline{c}_1^2-\overline{c}_2,\quad u_1^2u_2^2=\overline{c}_2,\quad
u_1u_2^3=\overline{c}_1^2-3\overline{c}_2,\quad u_2^4=5\overline{c}_2-\overline{c}_1^2,\cr
&u_1^3\cdot F=0,\quad u_1^2u_2\cdot F=-\overline{c}_1\cdot F,\quad u_1u_2^2\cdot F=0,
\quad u_2^3\cdot F=0,\cr}
$$ where $F$ is any divisor in $\Pic (X).$\noindent

Using this table, we obtain easily

$$(2u_1+u_2)^3\cdot Z = m_0(13\,\overline{c}_1^2-9\,\overline{c}_2) + 12\,t_0\,\overline{c}_1^2,$$ 
moreover we have 
${t_0/m_0} =\overline\theta_2$,
 hence
$$
(2u_1+u_2)^3\cdot Z =  m_0((13+12\,\overline\theta_2)\,\overline{c}_1^2 - 9\,\overline{c}_2)>0
.$$
As in [DEG00] Proposition 3.4., we conclude that
the restriction $\cO_{\overline{X}_2}(1)_{| Z}$ is big.  
Consequently,  by Theorem~1.1.3, 
 every nonconstant entire curve $f:\bC\to X$ is such 
that $f_{[2]}(\bC)$ is contained in the base locus of 
$\cO_{\overline{X}_2}(l)\otimes\pi^\star_{2,0} \cO (-A)_{| Z}$ for $l$ large. 
This base locus is at most $2$-dimensional,
and projects onto a proper algebraic subvariety $Y$ of~$\overline{X}_1$. 
Therefore $f_{[1]}(\bC)$ is contained in $Y$, and the Theorem
is proved. \finpr

\ssection 1.4.  Complement of curves in $\bP^2$ |In this section we will 
consider the case $ (\bP^2,C) $ where $C$
is a plane curve of degree $d$. We will estimate the associated $2$-jet 
log-threshold. We start with a vanishing theorem of
log-symmetric differentials (similar of that of Sakai [Sa78]).

\th 1.4.1. Lemma|
 Let $C$ be a smooth curve of degree $d$
in~$\bP^2$, $m$ a nonnegative integer and $k\in\bZ$. Then
\smallskip 
$H^0(\bP^2,S^m\overline{T}^\star_{\bP^2}\otimes\cO(k))=0$~~ for all
$k\le\min(m-1,d-2)$.
In particular, for $d\ge 4, (\bP^2,C) $ is of log-general type and we have
the estimate
 $$
(d-3)\overline\theta_{1,m}\ge{\min(1,(d-1)/m )}.$$

\finth

\dem 
We consider the natural ramified 
covering $X\subset\bP^3 $ over $\bP^2$ associated to $C$ (if $C$ is given by
$P(z_0,z_1,z_2)=0$, then $X$ is defined by $z_3^d=P(z_0,z_1,z_2)$), 
let $L$ be the hyperplane section in $X$ over $C$.
Then we have an injective morphism (by taking pullbacks)
$$H^0(\bP^2,S^m\Omega_{\bP^2} (\log C) (k))\hookrightarrow H^0(X,S^m\Omega_X (\log L)
(k)),$$ the last group is contained in $H^0(X,S^m\Omega_X(m+k))$ which vanishes
for $k\le\min(m-1,d-2)$ by Lemma 5.1 in [DEG00]. 

Now we have $\overline{K}_{\bP^2}=\cO_{\bP^2}(d-3)$, consequently, there
are no nonzero sections in $H^0(\bP^2,S^m\overline{T}^\star_{\bP^2}\otimes\cO(t\overline{K}_{\bP^2}))$ unless
$t(d-3)\ge \min(m,d-1)$,
whence the lower bound for $\overline\theta_{1,m}$. \finpr

Using the above vanishing lemma we obtain a lower bound on the 2-jet
log-threshold 

\th 1.4.2. Lemma|
Let $C$ is a curve of degree $d\geq 4$ in $\bP^2$. Suppose that
the 2-jet base locus $\overline{B}_2$ associated with $(X,D)=(\bP^2,C)$ is
of the form $Z_\sigma=Z_0\cup\Gamma_2,$ where $Z_0$ an irreducible section 
and $Z_\sigma$ is the set of zeros of a section $\sigma\in H^0(\overline{X}_2,\cO_{\overline{X}_2}(m_0)\otimes\cO(t_0\,\overline{K}_X)).$ 
Then for $m_0\geq 6$ we have the estimate
 $$\overline\theta_2=\overline\theta_{2,m_0}\geq\max ({-1\over m_0};\min
 ({1\over{2(d-3)}}-{1\over 6}, {{d-1}\over{2m_0(p_0-1)(d-3)}}-{1\over
   6}),$$
where $p_0=[{m_0\over 3}].$
\finth
\dem 
Observe that $\sigma$
can be considered as a global holomorphic section of the bundle
$E_{2,m_0}\overline{T}^\star_X\otimes \cO(t_0\overline{K}_X)$. 
By the filtration of $E_{2,m_0}\overline{T}^\star_X$, we have a short exact 
sequence
$$0\to S^m_0 \overline{T}^\star_X\to E_{2,m_0}\overline{T}^\star_X \to E_{2,m_0-3}\overline{T}^\star_X\otimes
\cO(\overline{K}_X)\to 0.$$ Multiply all terms by $\cO(t_0 \overline{K}_X)$ 
and consider the associated sequence in cohomology. As $t_0<0$ and by Lemma 
1.4.1, the first $H^0$ group vanishes and  we get an injection
$$H^0\big( X, E_{2,m_0}\overline{T}^\star_X\otimes \cO(t_0 \overline{K}_X)\big) \longhook
H^0\big( X, E_{2,m_0-3}\overline{T}^\star_X\otimes \cO((t_0+1)\overline{K}_X)\big).$$
By assumption on $\overline{B}_2$ we must have $t_0+1\geq 0$, this gives the
firt part of the estimate.
Let now $m_0=3p_0+q_0$, $0\le q_0\le 2$ a positive integer. Then there is 
a (nonlinear) discriminant mapping
$$
\Delta:E_{2,m_0}\overline{T}^\star_X\otimes\cO (t_0\overline{K}_X)\to
S^{(p_0-1)(3p_0+2q_0)}\overline{T}^\star_X\otimes\cO ( (p_0+2t_0)(p_0-1)
\overline{K}_X ).
$$

In fact, we write an element of $E_{2,m}\overline{T}^\star_X$
in the form
$$
P(f)=\sum_{0\le j\le p}a_j\cdot f^{\prime\,3(p-j)+q}\,W^j
$$
where the $a_j$ is viewed as an element of
$S^{3(p-j)+q}\overline{T}^\star_X\otimes \overline{K}_X^j$, and
$ W \in \overline{K}_X^{-1}.$ 
The discriminant $\Delta (P)$ is calculated by interpreting $P$ as a
polynomial in the indeterminate $W$.

Applying this to $\sigma$, we obtain 
$2(p_0-1)t+p_0(p_0-1)\ge(p_0-1)(3p_0+2q_0)\theta_{1,
  (p_0-1)(3p_0+2q_0)}$ and this implies
$$
{t\over m}\ge{3p_0+2q_0\over 2m_0}\theta_{1,
  (p_0-1)(3p_0+2q_0)}-{p_0\over 2m_0}.$$
Using Lemma 1.4.1 again we obtain the remainer part of the estimate.
\finpr

We now turn to the question of the existence of $2$-jet
differentials of small degree. Recall from the filtration of the bundle
of $2$-jet differentials that we have an exact sequence
$$
0 \lra S^m\overline{T}^\star_X \lra
E_{2,m}\overline{T}^\star_X \buildo\lra^\Phi E_{2,m-3}\overline{T}^\star_X \otimes\overline{K}_X \to 0.
$$
We have the following ``proportionality'' Lemma

\th 1.4.3. Lemma|Let $(X,D)$ be a nonsingular surface of
log-general type. Then, for all sections
$$
P_i\in H^0(X,E_{2,m_i}\overline{T}^\star_X\otimes\cO_X(t_i\,\overline{K}_X))
$$
with $m_i=3,\,4,\,5$ and $t_i\in\bQ$,  
the section $\beta_1 P_2-\beta_2P_1$ associated with $\beta_i=\Phi(P_i)$ 
can be considered as a section in
$$
H^0(X,S^{m_1+m_2-3}\overline{T}^\star_X\otimes\cO_X((1+t_1+t_2)\,\overline{K}_X)),
$$
and it vanishes when $1+t_1+t_2<(m_1+m_2-3)\overline{\theta}_{1,m_1+m_2-3}$.
\finth
\dem
The section $\beta_1 P_2-\beta_2P_1$ is contained in $H^0 (X,
E_{2,m_1+m_2-3}\overline{T}^\star_X
\otimes\cO_X((1+t_1+t_2)\,\overline{K}_X))$, its image by $\Phi$ is zero. Then
it can be considered as a section in  $H^0(X,S^{m_1+m_2-3}\overline{T}^\star_X\otimes\cO_X((1+t_1+t_2)\,\overline{K}_X))$.
\finpr

Now we have the following application of the proportionality Lemma
 
\th 1.4.4. Lemma|Let $C$ be a  generic curve of degree $d\ge 6$. Then
$$
\overline{\theta}_{2,m}\ge -{1\over 2m}+{1-(3+\epsilon)/2m \over d-3}\qquad\hbox{for $m=3,\,4,\,5,$}
$$ where $\epsilon:= d \mod 2$.

\finth

\dem
We consider the curve
$$
C_a=
\big\{z_0^{k_0}(z_0^{d-k_0}+a\,z_1^{k_1}z_2^{k_2})+z_1^d+z_2^d=0\big\},
$$
where $k_0,\,k_1,\,k_2,\ge 0$ are integers with
$\sum k_i=d$ and $a$ a complex number such that $C_a$ is non singular. 
We put
$$
s_0=z_0^{k_0}(z_0^{d-k_0}+a\,z_1^{k_1}z_2^{k_2}),\qquad 
s_i=z_i^d,\quad i=1,2.
$$
Using Nadel's method [Na89], we solve the linear system
$$
\sum_{0\le k\le 2} \wt\Gamma ^k_{ij} {\partial s_\ell\over \partial
z_k} = {\partial^2 s_\ell\over \partial z_i \partial z_j},\qquad 0\le
i,j,\ell \le 2,
$$
and get in this way a homogeneous meromorphic connection of degree 
$-1$ on $\bC^3$. One can check
that this connection descends to a partial projective meromorphic 
connection $\nabla=(\Gamma ^k_{ij})$ on $\bP^2$ such that $C_a$ is 
totally geodesic (see [EG96] and  [DEG97]), the pole divisor 
of the connection $\nabla$ is given by
$$
B=\big\{z_0 z_1 z_2 (d\,z_0^{k_1+k_2}+
ak_0z_1^{k_1}z_2^{k_2})=0\big\}.
$$ 
Then, this connection can be seen as a meromorphic connection on 
$\overline{T}_{\bP^2}$).
In fact, if we take two tangents vector fields $u$ and $v$ to $C_a$, the
vector field $ \nabla_uv$ is also tangent to $C_a$ by construction.

Consequently, by taking the Wronskian operator
$$
W_{\nabla}(f)=f'\wedge f''_{\nabla},\qquad f''=\nabla_{f'}f',
$$
 we get a section  $$P_1\in H^0 E_{2,3}\overline{T}^\star_{\bP^2} \otimes
 \cO(-\overline{K}_{\bP^2} )\otimes \cO(B)=
H^0 E_{2,3}\overline{T}^\star_{\bP^2} \otimes \cO(t_1\overline{K}_{\bP^2} ),$$ where
$t_1={{3+k_1+k_2}\over d-3}-1$. Remark that $p=3+k_1+k_2$ can be taken 
equal to any integer in $[3,d+3]$.  

We take $p=[{d+1 \over 2}]$ (the biggest integer less or equal to ${d+1 \over 2}$), 
so that
$$
{1\over 2}+t_1={3 +\varepsilon\over 2(d-3)},\quad
\hbox{where}\quad \varepsilon=d\mod 2,~~\varepsilon\in\{0,1\}.$$
The integer $p$ must be at least equal to $3$, thus our choice is 
permitted if $d\ge 6$. We claim that ${\bP^2}$ has no non 
trivial section in
$$
H^0({\bP^2},E_{2,m}\overline{T}^\star_{\bP^2}\otimes\cO(t\overline{K}_{\bP^2})), \qquad m=m_2=3,\,4,\,5
$$
if ${1\over 2}+t<{ m-3/ 2 -\varepsilon/2\over d-3}$. 
We assume the contrary, so that $P_2\in
H^0({\bP^2},E_{2,m}\overline{T}^\star_{\bP^2}\otimes\cO(t\overline{K}_{\bP^2})), P_2\neq0$. Then,
for $m_1=3$, $m_2=m$ and $t_2=t$, our choices imply
$$
1+t_1+t_2<{ m\over d-3}\le(m_1+m_2-3)\overline\theta_{1,m_1+m_2-3},
$$
as $\overline\theta_{1,m}\ge {1 \over d-3}$, for all
$m=3,\,4,\,5$ and $d\geq 6$ (by 1.4.1).
By Lemma~1.4.3, we get a
meromorphic connection with logarithmic pole on $C_a$ associated with a 
Wronskian operator
$P_2/\beta_2=P_1/\beta_1$. As $P_1/\beta_1$ is an irreducible
fraction with $\div\beta_1=B$, we conclude that
$\beta_2/\beta_1\in H^0({\bP^2},S^{m-3}\overline{T}^\star_{\bP^2}\otimes
\cO((t_2-t_1)\overline{K}_{\bP^2}))$ must be holomorphic, hence 
$$
t_2\ge t_1+(m-3)\overline\theta_{1,m-3}\ge t_1+{(m-3)\over d-3},$$ 
On the other hand
$$
t_2=t<-{1\over 2}+{ m-3/2-\varepsilon/2\over d-3}
=t_1+{ m-3-\varepsilon\over d-3},
$$
which yieldes a contradiction. By the Zariski semicontinuity of cohomology, 
the group
$$
H^0({\bP^2},E_{2,m}\overline{T}^\star_{\bP^2}\otimes\cO(t\overline{K}_{\bP^2}))
$$
vanishes for a  generic curve~$C$ of degree $d\geq 6$, unless
$$
{t\over m}\ge-{1\over 2m}+{1-(3+\varepsilon)/{2m}\over d-3}.
$$
this yields the estimate.\finpr

\finpr

As a corollary we obtain the main result of this first part

\th 1.4.5. Theorem|

Every non degenerate holomorphic entire map into the complement of a 
generic curve of degree $d\ge 15$ in $\bP^2$ is a leaf of a multi-foliation
on $\bP^2$.
\finth

\dem 
Let $C$ be a genenic curve of degree $d\geq 15$ in $\bP^2.$ 
If $\overline{B}_2$ is not as in Lemma 1.4.2, then we are done (we have
two independent sections), so
suppose that $\overline{B}_2=Z\cup \Gamma_2$ (with the notations of lemma
1.4.2). If $m_0=3,4,5,6,7$ (we have clearly
$m_0\geq 3$) we apply 1.4.4 and 1.4.2 to get the estimate
on the $2$-jet threshold
 $$\overline{\theta}_2 \geq {(3-\epsilon)\over{6(d-3)}}-{1\over 6},$$
where $\epsilon= d \mod 2$.
According to this estimate, 
$$(13+12\overline{\theta}_2)\overline{c}_1^2- 9\overline{c}_2 \geq (d-3)(2d-27-2\epsilon) -27$$
this is positive when $d\geq 15$ and we can apply Theorem 1.3.2
to obtain the statement.

When $m_0\geq 8$ we apply the estimate in lemma 1.4.2, we obtain 
$\overline{\theta}_2\geq -{1\over 8}.$  It is easy to verify that
$(13-3/2)\overline{c}_1^2- 9\overline{c}_2= (d-3)(2,5d-34,5)-27$ is positive for
$d\geq 15$ and we can again  apply Theorem 1.3.2 .

 \finpr

\section 2. Entire leaves of foliations on log-general type surfaces|
In this part we will generalize the main result in [Mc98], we will
follow basically the strategy in [Br99] with a little improvement due to
the "convenience" of the logarithmic formalism.

\ssection 2.1. Singularities of foliations on surfaces|Let $X$ be a compact 
complex surface. Recall from [GM87] that we have a 
bijective correspondance between a holomorphic foliation $\cF$
on the surface $X$ with isolated singularities and a locally free subsheaf 
of the tangent sheaf noted $T_\cF$. In this case we have an exact sequence 
of sheaves
$$0\lra T_\cF \lra T_X \lra N_\cF. I_Z \lra 0,$$
where $N_\cF$ is called the normal bundle of the foliation and $I_Z$ an ideal
supported on the singularity set $Z$ of $\cF$.

The elements of $Z$ are the points where local vector fields defining $\cF$ 
vanish. Suppose that $\cF$ is given
around a singular point $p$ by a vector field $v$, then we note by $\lambda_1$
and $\lambda_2$ the eigenvalues of the linear part of $v$ and we make the
following 

\th 2.1.1 Definition|
The singularity $p$ is called reduced if the linear part $(Dv)(p)$ is nonzero
(say $\lambda_2\neq 0$) and the quotient $ \lambda={\lambda_1\over\lambda_2 }$
is not a positive rational number.
\finth

A reduced singularity at $p$ is called nondegenerate if $\lambda_1\lambda_2\neq
0$  and a saddle-node otherwise.
The importance of those singularities comes from the following well-known
theorem

\th 2.1.2 Theorem([Se68])| 
There exist a sequence of blow-ups $\sigma : \tilde{X}\to X$ such that the 
foliation $ \sigma^\star \cF$ has only reduced singularities.\finth

Now let $D$ be a normal crossing divisor on $X$, we say that the foliation 
$\cF $ defines a logarithmic foliation on $(X,D)$ if  $\cF$ is tangent to each 
component of $D$. The sheaf 
injection from $T_\cF \hookrightarrow T_X$ factors to a sheaf injection
$$0 \lra T_\cF \hookrightarrow T_X(-\log D) \lra N_\cF(-D). I_{Z'}\lra 0,$$
where $Z'$ is the set of logarithmic singularities of $\cF$ which is obviously
contained in $Z$.   The bundle $N_\cF(-D)$ will be called the logarithmic
normal bundle of $\cF$, and denoted by $\overline{N}_\cF$.

To a logarithmic foliation $(\cF,D)$ we associate a (singular) surface 
$\tilde{X}$ in the $1$-jet logarithmic space $\overline{X}_1$, called the 
logarithmic graph of $(\cF,D)$, which consists of the adherence of the 
liftings of all leaves or equivalently the blow-up of $X$ along the ideal 
$I_{Z'}$.  

When $\cF$ has only reduced singularities, and if $p$ is  a nondegenerate 
singularity in $Z'$, then the surface $\tilde{X}$ is smooth around $p$ and 
isomorphic around $\pi^{-1}(p)$ to the blow up  of $X$ at $p$. In the case  
$p$ is a saddle node of multiplicity $d$ (i.e., the first terms of $v$ in 
local coordinates are $z{\partial\over\partial{z}} + 
w^d{\partial\over\partial{w}}),$ then $\tilde{X}$ has a singular point of
type $A_{d-1}$ (i.e., it is given in local coordinates by $z^d=xy$).

We finish this section by the following useful natural formula
$$\cO_{\overline{X}_1}(-1)|_{\tilde{X}} = \pi^\star(T_\cF)\otimes\cO(\sum_{p\in Z'} d_pE_p),$$
where $d_p$ is the multiplicity of $p$ and $E_p$ the fibre $\pi^{-1}(p)$.

\ssection 2.2. The log-tautological inequality and consequences|
Recall from [Mc98] that to a holomorphic curve $f : \bC\to X$, 
where $X$ is supposed to be endowed with a  K\"ahler form $\omega$, 
we can associate a closed positive current in the following way: 
for every 2-form $\eta$ and $r>0$ we define 
$$T_r(\eta)=\int_0^r{dt\over t}\int_{D(t)}f^\star\eta,$$
where $D(t)$ is the disc of radius $t$. Then we consider the positive currents
$\Phi_r$ defined by
$$\Phi_r(\eta):= {{T_r(\eta)}\over T_r(\omega)} \;  \forall \eta\in A^2(X).$$
The family $\{\Phi_r\}_{r>0}$ is bounded and we can see easily that there is
a closed positive current $\Phi$ in its adherence. 
When $f(\bC)$ is not contained in a hypersurface $Y$, we prove, using
the Lelong-Poincar\'e formula, that $\Phi$ has positive intersection with $Y$.
As a consequence, if $f$ is nondegenerated, then the current $\Phi$
is actually numerically effective.

Remark that this contruction is independent from the dimension of $X$.
Then,  
we can associate a positive current $\Phi_1$  on the $\overline{X}_1$ to the
curve $f_1$ (the lifting of  $f$). Now  if we suppose that $f$ intersects 
$D$ in finite set, we have the following logarithmic tautological
inequality (see [Vo99] and [Mc99]) 
 $$\cO_{\overline{X}_1} (-1). \Phi_1 \geq 0,$$
As a consequence of this inequlity we have $\pi_\star\Phi_1=\Phi.$ 
From now on we will suppose that there is a logarithmic foliation $\cF$ on 
$(X,D)$ with reduced singularities such that $D$ is the union of its 
algebraic leaves and $f$ is a Zariski dense entire leaf. 
Let $\nu(\Phi,p)= [\Phi_1].[d_pE_p],$ where $d_p$ and $E_p$ are defined in 
the previous section, then we have the following 

\th 2.2.1. Observation|
The logarithmic tautological inequality applied to the triple $(X,\cF,D)$ 
implies the refined tautological inequality.
\finth

\dem 
 We apply the logarithmic
 tautological inequality to $f_1$ which gives
$$
 \pi^\star(T_\cF)\otimes\cO(\sum_{p\in Z'} d_pE_p). \Phi_1 \geq 0.$$
As $\pi_\star\Phi_1=\Phi$, then we obtain

$$T_\cF . \Phi \geq -\sum_{p\in Z'} \nu(\Phi,p).$$

This last inequality is exactly the refined tautological inequality
because the intersection of algebraic leaves are not counted in $Z',$ 
in fact those points are smooth from the logarithmic point of view.
\finpr 

Now we iterate this construction: $(X,D)$ will be replaced by
$(\tilde{X},\tilde{D})$, where  $\tilde{D}$ is the union of all algebraic
leaves of the induced foliation on $\tilde{X}$ (we replace $\tilde{X}$
by its desingularisation if necessary). 
The last step is almost the same as in [Mc98] and [Br99]. 
Let $(X^{(n)}, \cF^{(n)})$ be the foliated
surface obtained after $n$ iterations of this process and let $\pi^{(n)}$ 
denotes the canonical morphism from  $X^{(n)}$ to $X^{(0)}=X$. 
As the singularities are reduced we have $$\cF^{(n)}=(\pi^{(n)})^\star(\cF).$$ 
Suppose (for simplicity) that we start from $Z'=\{p\}$ only one point of
multiplicity $d_p$, then 
$\cF^{(n)}$ have at most two logarithmic singularities: surely $q_1^{n}$ of 
multiplicity $d_p$ and  probably
$q_2^{n}$ of muliplicity $1$ so the log-tautological inequality gives

$$T_\cF.\Phi=T_{ \cF^{(n)}}.\Phi^{(n)}\geq -\nu(\Phi^{(n)},q_1^{(n)})-\nu(\Phi^{(n)},q_2^{(n)}).$$

Now we prove that $\nu(\Phi^{(n)},q_i^{(n)})$ tends to zero as $n$ tends to 
infinity, in fact, by comparing $\Phi^{(n)}$ and $\Phi$, we obtain the inequality
$$0\leq [\Phi^{(n)}]^2\leq [\Phi]^2-\sum_{j=0,i=1,2}^{n-1}\nu(\Phi^{(n)},q_i^{(n)}).$$
Consequently, we obtain
\th 2.2.2. Theorem([Mc98])|
Let $\cF$ a holomorphic foliation (with reduced singularities) on a 
compact surface $X$ and $\Phi$ the current associated to a Zariski dence
entire leaf. Then we have the intersection inequality

$$T_\cF . \Phi \geq 0.$$
\finth
 
\ssection 2.3. Positivity of the log-normal bundle on leaves|
Recall first the following result from [Br99]

\th 2.3.1. Theorem([Br99])|

Let $\cF$ be a holomorphic foliation (with reduced singularities) on a 
compact surface $X$ and $\Phi$ be an diffuse current (i.e., with zero Lelong 
numbers exept at a finite set of points. Suppose that $\Phi$ is 
$\cF$-invariant. Then, we have the intersection 
inequality  $$N_\cF(-D) . \Phi\geq 0, \hbox{ forall invariant divisor } D.$$   

\finth

The proof of the previous theorem consists of an explicit construction of a
closed $2$-form which represents the Chern class of $N_\cF(-D)$. The 
intersection is computed by integrating this form along the support
of $\Phi$ which is a union of leaves. This integration is concentrated 
around singularities for which structure and holonomy are well understood.

Now using this, the generalisation of the positivity of the normal bundle
on the current associated with a Zariski dense leaf to the log-case is
immediate from the following

\th 2.3.2. Corollary|

Let $(\cF,D)$ a holomorphic log-foliation with reduced singularities 
and $\Phi$ an $\cF$-invariant current such that the support of $\Phi_{\rm alg}$ is
contained in $D$. Then we have the following intersection inequality for the
logarithmic normal bundle
 
 $$ \overline{N}_\cF . \Phi\geq 0 .$$
\finth

\dem
 By Theorem 2.3.1. it remains to prove that $ \overline{N}_\cF . \Phi_{\rm
 alg} \geq 0 $. We use the same observation as in [Br99]: let $C$ be a 
component of $\Phi_{\rm alg}$, then  
$$ N_\cF(-D) . C= C.C + Z(C ,\cF) - D.C ,$$
where $ Z(C ,\cF)$ is the total multiplicity of the singularities of $\cF$
along $C$ (cf [Br97] lemme 3). This number is at least equal to the
intersectin $(D-C).C$, so we obtain  $ N_\cF(-D) . C \geq 0$. This is
True for evry component in the support of $ \Phi_{\rm alg}$, which
concludes the proof.
\finpr

As a consequence we have the following

\th 2.3.3.Theorem|
Let $(X,D)$ be a surface of log-general type with a logarithmic
foliation  $(\cF,D)$, then an entire leaf of $\cF$ must be degenerated.
\finth
\dem By the Seidenberg theorem we can suppose that $\cF$ has only reduced
singularities. Suppose that $\cF$ has a Zariski dence entire leaf.
Then, by theorem 2.2.2 and corollary 2.3.2,  as 
$K^{-1}_X=T_\cF \otimes N_\cF$, we obtain
 $$K_X\otimes\cO(D') . \Phi\leq 0,$$
where $D'$ is the union between $D$ and the support of $\Phi_{\rm alg}$.  
As $K_X\otimes\cO(D)\hookrightarrow K_X\otimes\cO(D'),$ the latter bundle
is big and, hence, has the decomposition $\cO(A+E)$, where $A$ is ample
and $E$ is effective. But $\Phi$ is numerically effective, so we obtain an 
obvious contradiction.
\finpr  

\ssection 2.4. Entire leaves on a surface of log general type|
In this section we will generalise the main theorem in [Mc98]. We must
consider now non necessarilyy logarithmic foliation.

\th 2.4.1. Lemma|
Let $\cF$ be a foliation  on a surface $X$ with a 
non-invariant curve $C$. Suppose that there is a Zariski-dense entire curve
 $f : \bC\to X\setminus C$. Then there exists a sequence of blow ups 
$\sigma : \tilde{X}\to X$ such that if we denote by $\tilde{C}$
the strict transform of $C$ and $\tilde{f}$ the lifting to $\tilde{X}$ of $f$, 
the support of an associated current $\tilde{\Phi}$ to $\tilde{f}$ is disjoint
from $C$. In particular,    $\tilde{\Phi} .\tilde{C} = 0$. \finth

\dem
By the Seidenberg theorem, using a sequence of blow ups $\tilde{X}\to X$, 
we can reduce the  singularities of $\cF$ and suppose that the leaves of 
the induced foliation $\tilde{\cF}$ are smooth, we can also suppose that 
these leaves are transverse to $\tilde{C}$ by blowing up the tangency points. 
Let $\cal{L}$ be the leaf containig the image of $\tilde{f}$, then as $f$ is 
Zariski-dense, $\cal{L}$ intersects $\tilde{C}$ on at most one point 
($\cal{L}$ is parametrised by $\bC$ or $\bC^\star$). Blowing up this point if 
it exists, we can suppose that $\cal{L}$ does not intersect $\tilde{C}$. 
We will prove that the topological closure of $\cal{L}$, which we denote by 
$K$, does not intersect $\tilde{C}$. 
Remark that $K$ is a union of leaves and let $p$ a point on 
$\tilde{C}\cup K$, so there is a leaf ${\cal{L}}_p$ in $K$ passing by $p$. 
Now $\cal{L}$ accumulates on ${\cal{L}}_p$ and this leaf is trasverse to 
$\tilde{C}$, as consequence in a neighbourhood of $p$ the number of 
intersection points of $\cal{L}$ with $C$ is infinite, this is a
contradiction. 
Finally, the support of $\tilde{\Phi}$ is contained
in this closure so does not intersect $\tilde{C}$. \finpr

As a consequence of Theorem 2.3.3 and the previous Lemma we obtain
the following generalisation of the main theorem of [Mc98]

\th 2.4.2. Theorem|
Let $X$ be a surface with a foliation $\cF$ and a divisor $D$ such that
$(X,D)$ is of log-general type. Then every entire curve $f : \bC \to X\setminus
D$ contained in a leaf of $\cF$ is degenerated. \finth

\dem Suppose that we have a Zariski-dense entire curve contained in a leaf
of $\cF$. We make a sequence of blow-ups to reduce the singularities of $\cF$
and to make the non-invariant components of $D$ verify the conclusion of 
Lemma 2.4.1. We will use the same notation on the obtained surface. 
Let $D_1\subset D$ be the union of the components of $D$ transverse  to $\cF$,
and $D_2$ the union of algebraic leaves, then, using Lemma 2.4.1,
we obtain
$$(K_X + D_1 + D_2).\Phi = (K_X+ D_2). \Phi.$$
As $\cF$ induce a logarithmic foliation on $ (X,D_2 )$, using the proof
of Theorem 2.3.3 we obtain $ (K_X+ D_2). \Phi \leq 0,$ which implies
 $$(K_X + D_1 + D_2).\Phi\leq 0.$$ Now the divisor $K_X + D_1 + D_2$ is big
and $\Phi$ is nef, so we have a contradiction.\finpr

As a consequence, we obtain the following

\th 2.4.3 Corollary|
Let $(X,D)$ be a log-surface of log-general type such that its logarithmic
Chern classes verify $\overline{c}_1^2 > \overline{c}_2$. Then every entire
curve $f : \bC \to {X\setminus D} $ is degenerated.\finth

\dem 
 We can suppose that $(X,D)$ is minimal. Then we apply Riemann-Roch to 
 symmetric powers of $\overline{T}_X$, the Euler characteristic is positive
 with our condition on Chern classes.
 Using Serre duality and nefness of $K_X+D$, the $h^2$ term is bounded 
 by the $h^0$ term. As a consequence we get 
 $\overline{B}_1\neq \overline{X_1}$,
 and we apply Theorem 1.1.3 to get a foliation $\cF$ on a (singular) surface 
 $\tilde{X}$ in $X_1$ which is ramified over $X$ (the foliation is defined by
 the intersection of $V_1$ and the tangent of the irreducible component 
 of $\overline{B}_1$ which
 contains the lifting of $f$ to $1$-jets) . Let $\tilde{D}$ be the 
 divisor on $\tilde{X}$ over $D$, then $(\tilde{X},\tilde{D})$ is of 
 log-general type
 and $f$ can be lifted in $\tilde{X}\setminus\tilde{D}$ as a leaf of the 
 foliation $\cF$.
\finpr  

\section 3. Proof of the Main Theorem|
By the results of Theorem 1.4.5, if $C$ is a generic plane curve of degree
$d\geq 15$, then there is a ramified cover $\tilde{X}\subset \overline{X}_1$ 
over $X=\bP_\bC^2$ with a foliation $\cF$ such that every entire curve $f$ in 
$\bP_\bC^2\setminus C$ is such that $f_{[2]}(\bC )$ is contained in 
$\tilde{X}$ as leaf of $\cF$. Morover, $f_{[2]}(\bC )$ is 
contained in $\tilde{X}\setminus \tilde{C}$ where $\tilde{C}$ is the (reduced)
divisor in $\tilde{X}$ over $C$.

A log model of $(\tilde{X},\tilde{C})$  is clealy of log-general type. 
In fact, recall that  given two 
log-manifolds $(X,D_X)$ and $(Y,D_Y),$ a holomorphic map $\psi : X \to Y$ 
such that $\psi^{-1}D_Y \subset D_X$ (in the geometric sens) is called a 
log-morphism.
Such a morphism induces (see [Ii77]) a vector bundle morphism
 $$\psi^\star : \psi^\star\overline{T}^\star_Y \to \overline{T}^\star_X.$$
If $\psi$ is birationnal, then this morphism is clearly injective.
Thus we have a natural injection of sheaves
$\psi^\star\overline{K}_Y\hookrightarrow \overline{K}_X $. 

Now, by Theorem 2.4.2 every entire curve $f$ in $X\setminus C$ has its lifting
in $\tilde{X}\setminus\tilde{C}$ degenratated and then itself has its image
contained in an algebraic plane curve. Now every algebraic curve in $\bP^2$
intersects a very generic curve of degree 
$d\geq5$ in at least $3$ point (see [Si95]) and then $f$ 
is constant and $\bP_\bC^2\setminus C$ is hyperbolic and hyperbolically 
embedded in $\bP_\bC^2$ (see [Gr77]) .

\titre References|

{\eightpoint
\parskip=4pt plus 1pt minus 1pt
\lettre M-De78|

\divers DB00|Berteloot F., Duval J|Sur l'hperbolicit\'e de certains
compl\'ementaires| Preprint 2000|

\article Bo79|Bogomolov F.A|Holomorphic tensors and vector bundles on 
projective varieties|Math.\ USSR Izvestija|13|1979|499-555|
 
\article Bro78|Brody R|Compact manifolds and hyperbolicity|Trans.\
Amer.\ Math.\ Soc.|235|1978|213--219|

\article Br97|Brunella M|Feuilletages holomorphes sur les surfaces
complexes compactes|Ann. Sci. ENS|30|1997|569--594|

\article Br99|Brunella M|Courbes enti\`eres et feuilletages
holomorphes|L'Enseignement Ma\-th\'e\-matique|45|1999|195--216|

\livre Br00|Brunella M|Birationnal geometry of foliation|First Latin
Amer. Cong. of Math, IMPA|2000|

\livre De95|Demailly J.-P|Algebraic criteria for Kobayashi
hyperbolic projective varieties and jet differentials|Proc. Sympos.
Pure Math., vol. 62, Amer. Math. Soc.,Providence,RI,1997|pp.285-360|

\livre DEG97|Demailly J.-P., El Goul J|Connexions m\'eromorphes
projectives et vari\'et\'es alg\'ebriques
hyperboliques|C.~R.\ Acad.\ Sci.\ Paris, t.324, S\'erie I, pp. 1385-1990 
|(1997)|
\article DEG00|Demailly J.-P., El Goul J|Hyperbolicity of generic surfaces of
high degree in projective 3-space|Amer.\ J.\ Math|122|2000|515--546|

\divers DL96|Dethloff G., Lu S|Logarithmic jet bundles and applications|
Preprint, to appear in Osaka J. Math. 2000| 

\article DSW92|Dethloff G., Schumacher G., Wong P.M|Hyperbolicity
of the complement of plane algebraic curves|Amer.\ J.\ Math|117|1995|573--599|

\article DSW94|Dethloff G., Schumacher G., Wong P.M|On the
hyperbolicity of the complements of curves in Algebraic surfaces:
the three component case|Duke Math.\ Math.|78|1995|193--212|

\article EG96|El Goul J|Algebraic families of smooth hyperbolic
surfaces of low degree in $\bP^3_\bC$|Manuscripta
Math.|90|1996|521--532|

\article Gr75|Green M|Some examples and counterexamples in value
distrubution theory|Compos. Math.|30|1975|317--322|

\article Gr77|Green M|The hyperbolicity of the complement of $2n+1$
hyperplanes in general position in ${\bf P}^n_{\bf C}$ and related
results|Proc. Amer. Math. Soc.|66|1977|109--113|

\article GG80|Green M., Griffiths P|Two applications of algebraic
geometry to entire holomorphic mappings|The Chen Symposium 1979,
Proc.\ Inter.\ Sympos.\ Berkeley, CA, 1979, Springer-Verlag|{\rm, New
York}|1980|41--74|

\article GM87|Gomez-Mont X| Universal families of foliations by curves|
Ast\'erisque |150-151|1987|109--129|
 
\livre Hi66|Hirzebruch F|Topological methods in Algebraic Geometry|
Grundl.\ Math.\ Wiss.{\bf 131}, Spriger, Heidelberg|1966|

\livre Ii77|Iitaka S| On the logarithmic kodaira dimension of algebraic varieties|
Complex Anal. and Alg. Geom. (ed. W.L. Baily, T. Shioda), Ianami Shoten,
175-189|1977| 

\livre Ko70|Kobayashi S|Hyperbolic manifolds and holomorphic
mappings|Marcel Dek\-ker, New York|1970|

\livre KR|Kobayashi R|Einstein-K\"ahler metric on open algebraic surface
of general type| Tohoku Math. J., 43-77|1985|

\article Lu91|Lu S|On meromorphic maps into varieties of log-general type|
Proceedings of Symposia in Pure Maths, Amer. Math. Soc|52(2)|1991|305--333|

\article McQ98|McQuillan M| Diophantine approximations and foliations|
Publ.\ Math.\ IHES|87|1998|121--174| 

\divers McQ99|McQuillan M| Non commutative Mori theory| Preprint|

\article MN94|Masuda, K., Noguchi, N|A construction of hyperbolic
hypersurfaces of $\bP^n$| Math. Ann. |304|1996|339--362| 

\article Na89|Nadel A.M|Hyperbolic surfaces in $\bP^3$| Duke Math. J.|
58|1989|749--771|

\article No77|Noguchi J|Meromorphic mappings into a compact complex space|
Hiroshima Math. J. |7|1977|411--425|

\livre No86|Noguchi, J|Logarithmic jet spaces and extensions of de
Franchis'Theorem|Contributions to Several Complex Variables (Conference in
Honor of W. Stoll, Notre Dame 1984), Aspects of Math., 227-249|1986|

\livre Sa78|Sakai F|Symmetric powers of the cotangent bundle and
classification of algebraic varieties|Proc.\ Copenhagen Meeting in
Alg. Geom.|1978|

\article Se68|Seidenberg A| Reduction of singularities of the differential
equation $AdY=BdX$|Amer.\ J.\ of Math.|90|1968|248--269|

\livre SY95|Siu Y.-T., Yeung S.K|Hyperbolicity of the complement of a
generic smooth curve of high degree in the complex projective plane| to
appear in Inventiones Math.|1996|

\livre TY86|Tian G., Yau S.T| Existence of K\"ahler-Einstein metrics
on complete K\"aler manifolds and their applications to algebraic 
geometry| in "Mathematical Aspects of String Theory" edited by S.T. Yau 
World Scientific|1986|

\divers Vo99|Vojta P|On the ABC conjecture and diophantine approximation by
rational points| preprint (1999)|

\article Zai87|Zaidenberg M|The complement of a generic
hypersurface of degree $2n$ in $\bC\bP^n$ is not hyperbolic|
Siberian Math.\ J.|28|1987|425--432|

\article Zai89|Zaidenberg M|Stability of hyperbolic embeddedness
and construction of examples|Math.\ USSR Sbornik|63|1989|351--361|
 
}

\parindent=0cm
\vskip20pt

Jawher El Goul\hfil\break
Universit\'e Paul Sabatier Toulouse III\hfil\break
D\'epartement de Math\'ematiques\hfil\break
118, route de Narbonne\hfil\break
31062 Toulouse, France\hfil\break
{\it e-mail:}\/ elgoul@picard.ups-tlse.fr

\immediate\closeout1
\end